\documentclass[12pt]{article}
\usepackage{amsfonts}
\usepackage{times}
\usepackage[left=1in,top=1in,right=1in]{geometry}
\usepackage{amsmath,amsthm,amssymb}
\usepackage{color}
\usepackage{latexsym}
\usepackage{cite}
\usepackage{hyperref}
\usepackage{latexsym}

\theoremstyle{plain}

\theoremstyle{plain}

\newtheorem{definition}{Definition}

\title{Comments on ``On Semi-Invariant Submanifolds of Almost Complex Contact Metric Manifolds"\\
	(Facta Universitatis, Series: Mathematics and Informatics  v. 31, no. 4, 851-862, 2016)}
\author{Aysel TURGUT VANLI  \\
	Dept. of Math. Fac. of Sci. Gazi University  TURKEI }
\date{}

\begin{document}

\maketitle

	\begin{abstract}
In this paper, we present the invalidities of the results in  \cite{cumalifeyzasemi}, because of their definition of a semi-invariant submanifold of an almost complex contact metric manifold is not true . 
	\end{abstract}
	\noindent\textbf{Keywords:} complex contact manifold, semi-invariant submanifold, normal complex contact metric manifold 
	
	\noindent \textbf{2010 AMS Mathematics Subject Classification:} 53C15,
	53C25, 53D10 
	\newline 
	\\
	
%
%
The definition of  a semi-invariant submanifold of an almost complex contact metric manifold is given as follows in \cite{cumalifeyzasemi}; 
\begin{definition}
	Let $ (\bar{M},\bar{X},\bar{Y},\bar{x},\bar{y},\bar{g},\bar{H}=\bar{G}J) $ be a an almost complex contact metric manifold . A submanifold $ M $ is a semi-invariant submanifold of almost complex contact metric manifold $ \bar{M} $, if there is $ (D,D^{\perp}) $ orthogonal distribution on $ M $ providing the 	following conditions;
	\begin{enumerate}
		\item $ TM=D\oplus D^{\perp} $
		\item $ D $ is invariant according to $ \bar{G} $ , that is, $\bar{G}D_z = D_z  $ for any $ z \in M$.
		\item $ D^{\perp}  $ is anti-invariant according to $\bar{G}$ , that is $\bar{G}D^{\perp}_z\subset T_zM^{\perp}$, for any $ z \in M, $
		where $ D $ and $  D^{\perp} $ distributions are horizontal and vertical distributions respectively.
	\end{enumerate}
	
\end{definition}
We have following remarks for this definition: 
\begin{itemize}
	\item It is not sufficient to give the definition only with the condition of $ \bar{G} $ tensor. Either conditions with $ \bar{H} $ or $ \bar{J} $ must be given.  
	\item A complex almost contact metric structure is defined on a complex manifold. For details we refer to reader \cite{BL2010}. But in this definition, it is not understand that the ambient manifold is complex or real. Also same case we have same situation for the submanifold $ M $.  
	\item The characteristic vector fields $ \bar{X}, \bar{Y} $ are not determined. They could be tangent or normal to the submanifold $ M $. In real case this two conditions are researched by different perspective and results could be different. We referee to reader Yano \cite{yanoCR} for details. But the authors don't give any details about the structure vector fields $\bar{X},\bar{Y} $, so there is a complexity.
	 \item
	 For any $ W $ vector fields tangent to submanifold $ M $ there is an expression as follow \cite{cumalifeyzasemi}:
	 \begin{equation} \label{w}
	 W=TW+RW+\bar{X}+\bar{Y} .
	 \end{equation} 
	 where $ T  $ and $ R $ are defined as the projection morphisms for  $ D $ and $  D^{\perp} $. \par 
We have two cases:\\ 
	\textbf{Case-1:} The structure vector fields $\bar{X},\bar{Y} $ are normal to $ M $. In this situation equation \ref{w} is not valid.
	 \\
		\textbf{Case-2:} The structure vector fields $\bar{X},\bar{Y} $ are tangent to $ M $. In this case $\bar{X}$ and $\bar{Y} $ could be in one of $ D , D^{\bot}$ or they could be orthogonal to these two distributions. 
		\subitem {i)} If $\bar{X},\bar{Y} \in \Gamma(D)$  then we have 
	\begin{eqnarray*} 
		\bar{g}(W,\bar{X})&=& \bar{g}(RW,\bar{X})+1\\
		\bar{g}(W,\bar{Y})&=& \bar{g}(RW,\bar{Y})+1
	\end{eqnarray*} 
	for all $ W \in \Gamma(TM) $. 
		\subitem {ii)}If $\bar{X},\bar{Y} \in \Gamma( D^{\bot})$ then we have
	\begin{eqnarray*} 
		\bar{g}(W,\bar{X})&=& \bar{g}(TW,\bar{X})+1\\
		\bar{g}(W,\bar{Y})&=& \bar{g}(TW,\bar{Y})+1
	\end{eqnarray*}
	for all $ W \in \Gamma (TM) $. 
	\subitem {iii)}	On the other hand if $X,Y \notin \Gamma(D\oplus D^{\bot})$  then we have
	\begin{eqnarray*}
		\bar{g}(W,\bar{X})&=& 1\\
		\bar{g}(W,\bar{Y})&=&1
	\end{eqnarray*}
	and thus we get $ W \in sp\{\bar{X},\bar{Y}\} $. Therefore we obtain $W=\bar{X}$ and $W=\bar{Y} $. So $ \bar{X}=\bar{Y} $, this is a  contradiction.\\
		
\end{itemize}
From all above remarks and determinations, we see that the definition is not true. On the other hand, all theorems in the paper are obtained from this definition. Thus, all results are not valid.

\end{document}